\documentclass[reqno, 11pt]{amsart}
\usepackage{amsmath,amssymb,supertabular}
\def\n{\noindent}
\def\hf{\hfil\break}

\begin{document}

\centerline {\bf
Axiomatic Definition of Limit of Real-valued Functions}
\bigskip

\n {\bf Abstract} We present a new way of organizing the few
mathematical statements which form introduction to Calculus: the
epsilon-delta characterization of the limit  is now  d e r i v e d
from four simple, intuitive and frequently used statements, which
we choose as axioms.
\bigskip

\n {\bf Introduction}

Several concepts of analysis, like real numbers and the Lebesgue measure,
have both axiomatic and non-axiomatic definitions. Apparently, until now
the limit concept had no axiomatic definition; in this note we present
such an axiomatic definition for the limit of  r e a l - v a l u e d
functions. The definition does not extend to limits in topological
spaces.  (The situation is somewhat different with metric spaces since a
function    $\varphi(x)$  with values in a metric space $(M,d)$ converges to a point
$y$ in $M$ if and only if the real-valued function  $f(x) = d(\varphi  (x), y)$
converges to zero).

One difference between the previous definitions and the new one is that
now the limit is defined as a mapping, previously as a relation (between
real-valued functions and real numbers). Because of that
we have now to define also the domain of the mapping,i.e. the class $F$ of
convergent functions. In this note we define $F$ axiomatically; it is
possible however to combine the two axioms that we use to define the limit
with a non-axiomatic definition of the class $F$ (as we did in [1]), or
to define $F$ simply as the largest class of functions on which a limit
exists (a limit satisfying our two axioms), but we think that the
axiomatic definition of $F$ provides more direct path into Calculus.

There are numerous kinds of limits in Calculus: limits of sequences of
real numbers; of real-valued functions of a real variable $x$ as $x$ tends to
$x_0 -$, to $x_0 +$, to $x_0$ , to $\infty  ,...;$ various limits of real-valued functions on
$R^n$, on a metric space; limits of Riemann and Riemann-Stieltjes sums. In
textbooks usually the limit is defined and the basic limit facts are listed for each
case separately; we shall in this note present
the axiomatic definition  for the left-hand limit of the real valued
function $f$ at the point $x_0$ and establish the equivalence with the standard
epsilon-delta definition; these statements are easily adapted to the
right-hand limit of $f$ at $x_0$ , or to the limit as $x$ tends to infinity , or
to minus infinity. However the adaptation  may not be obvious for the
two-sided limit, or for the limit of functions defined on $R^n$  or on a
metric space.  In the appendix we shall see that this adaptation is easy
and natural if the axioms are formulated for net convergence. However,
since the concept of nets is found mostly in advanced texts ([4],[7][10,[11]),
 (where it is used for the study of certain topological
spaces) and only rarely in introductory analysis texts ([2], [5], [6],[8], [9])
(where it is used to unify different kinds of limits
which appear in Calculus), we did not want to assume that the reader is
familiar with net convergence and therefore nets appear only in the
appendix, which, we hope, would be understandable even to a reader who
meets here the net convergence for the first time.

\n {\bf Axioms and their immediate consequences}

The two axioms defining the mapping ``left-hand limit at $x_0$" (denoted  $\lim\limits_{x \to x_0 -}$)
and the two axioms defining the domain of the mapping, i.e. defining
real-valued functions that ``converge as $x  \to x_0 -$", are interconnected.

The limit axioms are

\n (1) (Constants axiom) If $f(x) = c$ for all $x < x_0 $, then $\lim\limits_{x \to x_0 -} f(x) = c$.

\n (2) (Inequality axiom) If the functions $f$ and $g$ converge as $x  \to x_0-$  and if
$\lim\limits_{x\to x_0-} f(x) < \lim\limits_{x\to x_0-} g(x)$, then there exists
$ a < x_0$ such that  $f(x) < g(x)$  for $ a< x < x_0$ .

The inequality axiom is crucial, it plays here the role that the
functional equation of the exponential function has in elementary
analysis, or the role that countable additivity has in measure theory. By
raising this well-known statement to the piedestal of an axiom, we did
uncouple the traditional epsilon and delta.

Immediate consequences of the limit axioms:

By contradiction from (2) we derive

\n (3) (Inequality theorem) If the functions $f$ and $g$ converge as $x \to x_0- $ and
if there exists $a < x_0$ such that
$f(x) \le g(x)$  for $a < x < x_0$   then $\lim\limits_{x \to x_0 -} f(x) \le \lim\limits_{x \to x_0 -} g(x)$.

From (1) and (2) we obtain

\n (4)(i)  If $\lim\limits_{x \to x_0 -} f(x) = l$ and  $l' < l$ then there exists $a' < x_0$ such that $ l'
< f(x)$  for $a' \le x < x_0 $.

\n   (ii) If $\lim\limits_{x \to x_0 -} f(x) = l$ and  $l < l$" then there exists $a" < x$ such that
$f(x) < l"$  for $a" \le x < x_0 $.

From (4) it follows obviously

\n (5) If $f$ converges as $x \to x_0 -$, then there exists $a < x_0$ such that  $f$ is
bounded on the interval $[a, x_0 )$, and

\n (6) If     $\epsilon > 0$ and $\lim\limits_{x \to x_0 -} f(x) = l$, then there exists $a < x_0$ such that  $|f(x)
- l| < \epsilon$   for $a < x < x_0$ .

Finally, we can use (4) to prove the uniqueness of limit in the following
sense:

\n (7) If $\mathop{\lim'}\limits_{x \to x_0 -} $ and $\mathop{\lim^{\prime\prime}}\limits_{x \to x_0 -}$ are two mappings satisfying axioms (1) and (2),
having possibly different domains, and if $f$ belongs to both domains, then
$\mathop{\lim'}\limits_{x \to x_0 -}f = \mathop{\lim^{\prime\prime}}\limits_{x \to x_0 -}f$.

The uniqueness in stronger sense (the fact that also the domains of the
two mappings coincide) will follow from the epsilon-delta theorem.

The convergence axioms are

\n (8) (MB axiom). If there exists $a < x_0$ such that $f$ is monotone and bounded
on the interval $[a, x_0)$, then $f(x)$ converges as $x \to  x_0 -$

\n (9) (Sandwich axiom) Let  $a < x_0$  and let the functions $f, g$ and
$h$ satisfy $f(x) \le g(x) \le h(x)$ on the interval $[a, x_0 )$. If $f $ and $h$ converge
as $x \to x_0-$  and if $\lim\limits_{x\to x_0-} f(x) = \lim\limits_{x\to x_0-} h(x)$, then $g$ converges as $x  \to x_0- $.

An immediate consequence of (9) and (3) is

\n (10) If there exists $a < x_0$ such that $f(x) = g(x)$ on the interval $[a, x_0 )$,
then if $f$ converges as $x  \to x_0 $, so does $g$,
and $\lim\limits_{x\to x_0-} f(x) = \lim\limits_{x\to x_0-} g(x)$.

Keeping unchanged the two limit axioms and slightly modifying the two
convergence axioms we can define also the extended limit (limit with
values in the extended real line $R \cup\{\infty\}  \cup\{-\infty\}$) and the real-valued
functions that converge in the extended sense.

\n {\bf Equivalence with the epsilon-delta definition}

Consider the statements:

\n(i')$ f(x)$ converges as $x \to x_0 -$

\n (i") $\lim\limits_{x\to x_0-} f(x) = l$

\n (ii) for every  $\epsilon > 0$  there exists   $\delta =\delta  (\epsilon)$ positive such that
$|f(x) - l| \le \epsilon$   for  $x_0-\delta \le  x < x_0$ .

Epsilon-delta theorem. Both (i') and (i") hold if and only if (ii) holds.

\n  Proof.

If (i') and (i") hold, then so does (ii).
This is statement (6) which was already proved.
The converse is proved in two steps. First we show

\n (11) If both (i') and (ii) hold, then (i") holds.

By (ii) we have

\n (12)   $l -\epsilon   \le  f(x) \le l + \epsilon$    on the interval  $[x_0  -\delta , x_0)$
Since $f(x)$ converges as $x\to x_0 -$, we deduce from (12), (1) and (3) that
$l - \epsilon  \le  \lim\limits_{x\to x_0-}
 f(x) \le l + \epsilon$ .

\n  Since this holds for every positive $\epsilon$  we conclude
that $\lim\limits_{x\to x_0-} f(x) = l$, which proves (11). The
final step is to show

\n (13) (ii) implies (i').

Let

\n (14) $M(x) =\sup \{f(t)| x \le t < x_0 \}, m(x) = \inf \{f(t)| x \le t < x_0 \}$.

Obviously both $M $and $m$ are monotone functions, and -by (12)- they are both
bounded on some interval $[a, x_0)$. By the MB axiom they both converge as
$x \to  x_0- $, in other words both satisfy (i'). It follows from (12) that they
both satisfy (ii). So, by (11) both $M$ and $m$ satisfy (i") and we have

\n (15)  $\lim\limits_{x\to x_0-} M(x) = l$ and $\lim\limits_{x\to
x_0-} m(x) = l$.
\n Observing that the definition (14) implies

\n (16) $m(x) \le f(x) \le M(x)$ for $x < x_0 $,

we deduce from (15), (16) and the sandwich axiom that $f(x)$ converges as
$x \to x_0-$ , which ends the proof.

Once the epsilon-delta theorem is established, it can be applied to obtain
some other results as, for example,

\n (17) Theorem on convergence of the sum, product and quotient of two
convergent functions.

(As in some older textbooks, the epsilon-delta form of the
definition may be used to prove (i) if the functions $z'$ and $z"$
converge to zero, then so does  $z' + z"$, (ii) the function $L +
z(x)$ converges to $L$ if and only if the function $z$ converges
to zero, and (iii) if the function $b$ is bounded, and the
function $z$ converges to zero, then the product   $b z$ converges
to zero. From these three facts, by simple algebraic manipulations
one obtains  (17)).

\n {\bf Presentation to beginners}

To build a course for beginners in which the limit would be introduced
axiomatically one would need now to provide examples and exercises, to add
some details and to decide which statements and proofs given above to
place in some kind of appendix.

Obviously there are many ways such a course can be built, but whatever way
it is done, we think that it would be useful  b e f o r e  stating a
definition of the limit to  explicitly  introduce the relation
``$f(x)$ is ultimately less than $g(x)$ as $x$ tends to $...$";  perhaps to even
introduce a provisional notation for that relation, like  ``$f -\!\!\!< g$  at $x_0 -$"
(at $x_0 +$, at  $x_0$, at   $\infty, ...)$

\n {\bf Appendix}

There are two equivalent ways to unify different kinds of limits for real-
valued functions - one can introduce a ``filter" (or a ``filter-base") in
the domain $X$ of the function $f$ ([3],[10], [12]) or introduce a  ``direction" on
$X$ ([2], [5], [6], [8],[9],[10]).

The second approach is better suited for our axiomatic definition.
A ``direction" on a set $X= \{x,y,z,...\}$ is a binary relation on $X$ which
is transitive $(x \prec y $ and  $y  \prec z$  imply  $x \prec z)$, reflexive $(x\prec x$ for
every $x)$ and has the property  that for any two elements $x$ and $y$ in $X$
there exists an element $z$ in $X$ such that  $x  \prec z$ and $y\prec z$.  A directed set
$(X, \prec)$ is a set $X$ with a direction; if $a  \in X$, then the set  $\{x | a \prec x\}$
is called a tail of the directed set $(X, \prec)$; a real-valued function defined
on a tail of a directed set is called a real-valued net on that directed
set.  A real-valued net $f$ is monotone increasing (decreasing) if for some
tail $T$  and all $x,y$ in $T$ such that  $x  \prec y$   we have $f(x) \le f(y)
(f(x) \ge f(y))$; it is monotone if it is monotone increasing or monotone
decreasing; net $f$ is bounded means that there exists a tail on which the
function $f$ is bounded. The net $f$ converges to the real number $l$ if for
every  $\epsilon$  positive there exists a tail $T= T(\epsilon)$ such that if $x$ is in the tail $T$ then
$|f(x) - l| <  \epsilon$. This is the classical definition.

To formulate our axiomatic definition we introduce first a strict partial
order relation  $-\!\!\!<$  between real-valued nets on the directed set $(X,  \prec )$
defined by
$f -\!\!\!< g$  if and only if there exists a tail $T$ such that for all $x$ in $T$ we
have  $f(x) < g(x)$ (similarly $f-\!\!\!\le g$ corresponds to $f(x) \le g(x)$).

The limit and convergence axioms for nets on $(X, \prec )$ are

\n (1') (Constants axiom) If $f(x) = c$ for all $x$, then $\lim f = c$.

\n (2') (Inequality axiom) If the nets $f$ and $g$ converge and if
$\lim f < \lim g$, then  $f -\!\!\!< g$.

\n (8') (MB axiom).  Monotone bounded nets are convergent.

\n (9') (Sandwich axiom) Let the nets $f, g$ and $h$ satisfy
$f -\!\!\!\le g -\!\!\!\le h$. If $f$ and $h$ converge and if $\lim f = \lim h$,
then $g$ also converges.

On a given set $X$ a statement of the form ``as $x$ tends to ..."
denotes a direction. For example, if $(X,d)$ is a metric space,
and $ x_0  \in  X$, then  ``as $x$ tends to $x_0$" denotes the
direction defined  on the set  $X\backslash  \{ x_0 \}$ by   $x
\prec y$  if and only if    $d(y,x_0) \le d(x,x_0 )$. A
real-valued function $f$ on the set  $X\backslash  \{ x_0 \}$ is a
net on the directed set $(X\backslash \{x_0\}, \prec )$ . It is
easy to see that the  n e t  $f$ is monotone if and only if the f
u n c t i o n  $f$ has the property ``there exists a positive
number h and a monotone function m on the interval $(0,h]$ such
that $f(x) = m (d(x,x_0 ))$ for $x$ such that $d(x,x_0 ) < h$."
So, if we want to adapt the MB axiom to convergence on metric
spaces we have the choice of using nets or sacrificing simplicity
of the axiom. This holds even in the case when the metric space is
the real line, namely when directly defining the two-sided limit
at $x_0$; obviously the simplicity can in that case  be preserved
without introducing nets, just defining the two-sided limit with
the aid of the left and right limit.

\vfill\eject

\bigskip

\n Bogdan Baishanski\hf
Department of Mathematics\hf
Ohio State University\hf
231 West 18 Avenue\hf
Columbus, OHIO 43210
\hf
bogdan@math.ohio-state.edu
\end{document}